\documentclass[a4paper,12pt]{article}
\usepackage[english]{babel}
\usepackage[T2A]{fontenc}
\usepackage[cp1251]{inputenc}
\usepackage[english]{babel}
\usepackage{amsthm}
\usepackage[tbtags]{amsmath}
\usepackage{amsfonts,amssymb}
\begin{document}

\newtheorem{theorem}{Theorem}
\newtheorem{lemma}{Lemma}
\newtheorem{proposition}{Proposition}
\newtheorem{Cor}{Corollary}

\begin{center}
{\large\bf On Centrally Essential Subrings

 of Formal Triangular Matrix Rings}
\end{center}

\begin{center}
Oleg Lyubimtsev\footnote{Nizhny Novgorod State University, Department of Algebra, Geometry and Discrete Mathematics, scientific and educational mathematical center "Mathematics of technologies of the future", Nizhny Novgorod, Russia; email: oleg\_lyubimcev@mail.ru , https://orcid.org/0000-0002-8773-1943.},
Askar Tuganbaev\footnote{National Research University MPEI, Department of Higher Mathematics, Moscow, Russia; Lomonosov Moscow State University, Department of Higher Algebra, Moscow, Russia; tuganbaev@gmail.com , https://orcid.org/0000-0001-9582-3314.}
\end{center}
\textbf{Abstract.} A ring $R$ with a non-zero identity element is said to be centrally essential if for any non-zero element $a\in R$, there exist non-zero central elements $x,y\in R$ such that $ax = y$. We describe centrally essential rings in a large subclass of formal triangular matrix rings and in a subclass of the matrix ring $M_3(R)$ over the ring $R$.

\textbf{Keywords:} centrally essential ring, formal triangular matrix ring, essential submodule.

The work of Oleg Lyubimtsev is supported by Ministry of Education and Science of the Russian Federation, project H-044-0\_2026-2028. The study of Askar Tuganbaev is supported by Russian Science Foundation, project 26-11-00001.

\textbf{MSC2020 database 16D25, 16R99}

\section{Introduction}\label{sec1}

We consider only associative rings with a non-zero identity element and unitary modules. A ring $R$ is said to be \textbf{centrally essential} if for any non-zero element $a\in R$, there exist  non-zero central elements $x,y\in R$ such that $ax = y$. Centrally essential rings were studied in \cite{LT25IJ}, \cite{MT19b} and other papers. 

Let $R$, $S$ be  two rings  and let $M$ be a non-zero $R$-$S$-bimodule. A \textbf{formal triangular matrix ring} $T(R, S, M) = \begin{pmatrix}
R&M\\
0&S\\
\end{pmatrix}$ consists of formal matrices $\begin{pmatrix}
a&m\\
0&b\\
\end{pmatrix}$, where $a\in R$, $b\in S$ and $m\in M$ with component-wise addition and the following multiplication: 
$$
\begin{pmatrix}
a&m\\
0&b\\
\end{pmatrix}\begin{pmatrix}
a'&m'\\
0&b'\\
\end{pmatrix} = \begin{pmatrix}
aa'&am' + mb'\\
0&bb'\\
\end{pmatrix}.
$$
Formal triangular matrix rings are of great importance in ring theory and have been studied, for example, in \cite{BPR12}, \cite{HV99}. Let $R = S$ and denote the resulting ring as $T(R, M)$.
The ring $T(R, M)$ is not centrally essential, since it contains non-central idempotents; cf. \cite[Proposition 2.4]{MT19b}. We consider unital subrings of $T(R, M)$ with $M\neq 0$.

{\textbf{1.1. Remark.} If $K$ is a centrally essential subring of $T(R, M)$ and $Z(K)$ is the center
of $K$, then 
$Z(M) = \left\lbrace m\in M \mid A = \begin{pmatrix}
a&m\\
0&b\\
\end{pmatrix}\in Z(K) \right\rbrace\neq 0$. Otherwise, for the non-zero matrix 
$\begin{pmatrix}
0&m\\
0&0\\
\end{pmatrix}\in K$, we have $\begin{pmatrix}
0&m\\
0&0\\
\end{pmatrix}\cdot Z(K)\cap Z(K) = 0$; this is a contradiction, since the ring $K$ is centrally essential.

For an $R$-$R$-bimodule $M$, we set $\text{Ann}_MI = \{m\in M \,|\, mI = 0\}$, $\text{Ann}_IM = \{i\in I \,|\, iM = 0\}$, $C_R(M) = \{m\in M \,|\, mr = rm, \forall r\in R\}$. We recall that a submodule $A$ of a module $M$ is said to be \textbf{essential} (notation: $A\le_eM$) if $A\cap B\neq 0$ for any non-zero submodule $B$ of $M$.

The main result of the paper is Theorem 1.2.

{\textbf{1.2. Theorem.} Let $R$ be a ring, $Z(R) = Z$ be the center of $R$, and let $M$ be a non-zero $R$-$R$-bimodule. A unital subring $K$ of $T(R, M)$ is centrally essential if and only if 
$$
K = \left\lbrace\left.\begin{pmatrix}
a&m\\
0&b\\
\end{pmatrix}\right | \; a, b\in R;\, m\in M \,\, \mbox{and} \,\, a - b\in I\vartriangleleft R \right\rbrace,
$$
where $I = \left\lbrace a - b\in R \mid A = \begin{pmatrix}
a&m\\
0&b\\
\end{pmatrix}\in K \,\, \mbox{and} \,\, Z(M)(a - b) = 0 \right\rbrace$, and the subbimodule $N = \text{Ann}_MI\cap C_R(M)$ of $M$ and the ideal $J = \text{Ann}_RN$ of $R$ satisfy the following conditions:
\begin{enumerate}
\item[\textbf{1)}]
$\text{Ann}_IM\cap Z\le_e I_Z$ and $J\cap Z\le_e J_Z$;
\item[\textbf{2)}]
$N\le_eM_Z$. 
\end{enumerate}

\textbf{1.3. Corollary.} Let $R$ be a commutative ring. The ring $K$ is centrally essential if and only if $\text{Ann}_IM\le_eI$ and $\text{Ann}_MI\le_eM$. 

\textbf{1.4. Corollary.} Let $R$ be a commutative domain. The ring $K$ is centrally essential if and only if $\text{Ann}_IM\neq 0$ and $\text{Ann}_MI\le_eM$.

\textbf{1.5. Corollary.} If $R$ is a field and $K$ is a centrally essential ring, then $K$ is commutative.

For a ring $R$ and an $(R,R)$-bimodule $M$, the \textbf{trivial extension} $\Lambda = R\ltimes M$ of $R$ by $M$ (for example, see \cite{S86}) is the Cartesian product $R\times M$ with  componentwise addition and multiplication given by $(r,m)(r', m')=(rr', mr' + rm')$. Then $\Lambda\cong \left\lbrace\left.\begin{pmatrix}
r&m\\
0&r\\
\end{pmatrix}\right|\;\; r\in R, m\in M \right\rbrace$. As a formal triangular matrix subring defined above, we have $I = 0$ and $N = C_R(M)$.

\textbf{1.6. Corollary.} The ring $\Lambda$ is centrally essential if and only if $J\cap Z\le_e J_Z$ and $C_R(M)\le_e M_Z$.

We also consider a class of subrings in the full matrix ring $M_3(R)$ over the ring $R$. Let $s, t\in Z(R)$. We define rings
$$
H_{(s,t)}(R) = \left\lbrace\left.\begin{pmatrix}
a&0&0\\
c&d&e\\
0&0&f
\end{pmatrix}\right.\in M_3(R)\; \mid \; a, c, d, e, f\in R;\, a = d + sc, d = f + te \right\rbrace.
$$
The rings $H_{(s,t)}(R)$ were previously investigated for being $CU$  rings and $C\Delta$ rings; for more details see \cite{KHHC19}, \cite{DJHM25}.

\textbf{1.7. Theorem.} The ring $H_{(s,t)}(R)$ is centrally essential if and only if $R$ is a centrally essential ring.

\section{Proof of Theorem 1.2}

We recall that $Z$ denotes the center of the ring $R$, $Z(K)$ is the center of the ring $K$, $N = \text{Ann}_MI\cap C_R(M)$.

\textbf{2.1. Lemma.} 
$$
Z(K) = \left\lbrace\left.\begin{pmatrix}
a&m\\
0&b\\
\end{pmatrix}\right | \; a, b\in Z, \, a - b\in \text{Ann}_IM, \, m\in N \right\rbrace.
$$

\textbf{Proof.} If $A = \begin{pmatrix}
a&m\\
0&b\\
\end{pmatrix}\in Z(K)$, then 
$$
AA' = \begin{pmatrix}
aa'&am'+mb'\\
0&bb'\\
\end{pmatrix} = \begin{pmatrix}
a'a&a'm+m'b\\
0&b'b\\
\end{pmatrix} = A'A
$$
for any matrix $A' = \begin{pmatrix}
a'&m'\\
0&b'\\
\end{pmatrix}\in K$. Then $aa' = a'a$, $bb' = b'b$ for all $a', b'\in R$, i.e. $a, b\in Z$. In addition, we have $am' + mb' = a'm + m'b$. Since $a, b\in Z$, we have $(a - b)m' = a'm - mb'$. If $a' = b' = 1$, then $(a - b)m' = 0$ for any $m'\in M$, i.e. $a - b\in \text{Ann}_IM$. Then
$mb' = a'm$ for all $a', b'\in R$. If $a' = b'$, then $m\in C_R(M)$. Now, if $a' \neq b'$, then $m(a' - b') = 0$ and $m\in \text{Ann}_MI$. Therefore $m\in N$.

Conversely, let $A = \begin{pmatrix}
a&m\\
0&b\\
\end{pmatrix}\in K$ where $a, b\in Z$, $a - b\in \text{Ann}_IM$ and $m\in N$. For any matrix $A'\in K$, we have $aa' = a'a$, $bb' = b'b$ for all $a', b'\in R$, since $a, b\in Z$. Furthemore, it follows from $m\in N$ that $m(a' - b') = 0$. Finally, it follows from the conditions $a, b\in Z$ and $a - b\in \text{Ann}_IM$ that $m'b = am'$ Therefore, $AA' = A'A$ and $A\in Z(K)$.~$\square$

\textbf{2.2. Completion of the Proof of Theorem 1.2.}

Assume that $K$ is a centrally essential ring. 
It follows from Lemma 2.1 that for arbitrary matrices $A = \begin{pmatrix}
a&m\\
0&b\\
\end{pmatrix}\in Z(K)$ and $A' = \begin{pmatrix}
a'&m'\\
0&b'\\
\end{pmatrix}\in K$, we have $m(a' - b') = 0$. With the use of Remark 1.1, we have that the set
$I = \left\lbrace a - b\in R \mid A = \begin{pmatrix}
a&m\\
0&b\\
\end{pmatrix}\in K \,\, \mbox{and} \,\, Z(M)(a - b) = 0 \right\rbrace$ forms a proper ideal of $R$. It is easy to see that the set of matrices $K$ forms a ring. 

Now we verify that conditions 1) and 2) hold.

\textbf{1)} Let $0\neq i\in I\backslash (\text{Ann}_IM\cap Z)$. Then the matrix $
\begin{pmatrix}
i&0\\
0&0\\
\end{pmatrix}$ is not contained in $Z(K)$. Since $K$ is a centrally essential ring, we have
$$
0\neq \begin{pmatrix}
i&0\\
0&0\\
\end{pmatrix}\begin{pmatrix}
a&n\\
0&b\\
\end{pmatrix} = \begin{pmatrix}
ia&in\\
0&0\\
\end{pmatrix}
$$
for some central matrices $\begin{pmatrix}
a&n\\
0&b\\
\end{pmatrix}$ and $\begin{pmatrix}
ia&in\\
0&0\\
\end{pmatrix}$. It follows from Lemma 2.1 that $in = 0$. Then $0\neq ia\in \text{Ann}_IM\cap Z$, where $a\in Z$. 

Let $0\neq j\in J\backslash (J\cap Z)$. Then $\begin{pmatrix}
j&0\\
0&j\\
\end{pmatrix}\notin Z(K)$ and 
$$
0\neq \begin{pmatrix}
j&0\\
0&j\\
\end{pmatrix}\begin{pmatrix}
a&n\\
0&b\\
\end{pmatrix} = \begin{pmatrix}
ja&jn\\
0&jb\\
\end{pmatrix} \in Z(K),
$$
where $\begin{pmatrix}
a&n\\
0&b\\
\end{pmatrix}\in Z(K)$. Since $n\in N$, we have $jn = 0$. Therefore, $0\neq ja\in Z$ or $0\neq jb\in Z$ for $a, b\in Z$. 

\textbf{2)} Let $m\notin N$ for some $0\neq m\in M$. Then $\begin{pmatrix}
0&m\\
0&0\\
\end{pmatrix}\notin Z(K)$ and 
$$
0\neq \begin{pmatrix}
0&m\\
0&0\\
\end{pmatrix}\begin{pmatrix}
a&n\\
0&b\\
\end{pmatrix} = \begin{pmatrix}
0&mb\\
0&0\\
\end{pmatrix} \in Z(K),
$$
where $0\neq mb\in N\cap mZ$.

Conversely, assume that the ring $K$ is of the specified form and conditions 1) and 2) hold. We prove that the ring $K$ is centrally essential.

Let $\begin{pmatrix}
a&m\\
0&b\\
\end{pmatrix}\notin Z(K)$. If $a = b = 0$, then it follows from \textbf{2)} that $0\neq mz\in N$ for some $z\in Z$. Then
$$
0\neq \begin{pmatrix}
0&m\\
0&0\\
\end{pmatrix}\begin{pmatrix}
z&0\\
0&z\\
\end{pmatrix} = \begin{pmatrix}
0&mz\\
0&0\\
\end{pmatrix}\in Z(K).
$$ 
Let $a\neq 0$. If $a\notin J$, then $am'\neq 0$ for some $m'\in N$. Then
$$
0\neq \begin{pmatrix}
a&m\\
0&b\\
\end{pmatrix}\begin{pmatrix}
0&m'\\
0&0\\
\end{pmatrix} = \begin{pmatrix}
0&am'\\
0&0\\
\end{pmatrix},
$$ 
where $\begin{pmatrix}
0&m'\\
0&0\\
\end{pmatrix}\in Z(K)$. If $am'\notin N$, then it follows from \textbf{2)} that $0\neq (am')z\in N$ for some $z\in Z$. Then
$$
0\neq \begin{pmatrix}
a&m\\
0&b\\
\end{pmatrix}\begin{pmatrix}
0&m'z\\
0&0\\
\end{pmatrix} = \begin{pmatrix}
0&am'z\\
0&0\\
\end{pmatrix}\in Z(K).
$$
Let $a\in J$. It follows from \textbf{1)} that $0\neq ac\in J\cap Z$ for some $c\in Z$. We assume that $ac = bc$. If $mc\in N$, then 
$$
0\neq \begin{pmatrix}
a&m\\
0&b\\
\end{pmatrix}\begin{pmatrix}
c&0\\
0&c\\
\end{pmatrix} = \begin{pmatrix}
ac&mc\\
0&ac\\
\end{pmatrix}\in Z(K).
$$ 
Otherwise, it follows from \textbf{2)} that there exists $z\in Z$ such that $0\neq (mc)z\in N$. Then
$$
0\neq \begin{pmatrix}
a&m\\
0&b\\
\end{pmatrix}\begin{pmatrix}
cz&0\\
0&cz\\
\end{pmatrix} = \begin{pmatrix}
acz&mcz\\
0&acz\\
\end{pmatrix}\in Z(K).
$$ 
Let $ac - bc\neq 0$. It follows from \textbf{1)} that there exists $z\in Z$ such that $0\neq (ac - bc)z\in \text{Ann}_IM\cap Z$. Then $bcz\in Z$. If $mcz\in N$, then 
$$
0\neq \begin{pmatrix}
a&m\\
0&b\\
\end{pmatrix}\begin{pmatrix}
cz&0\\
0&cz\\
\end{pmatrix} = \begin{pmatrix}
acz&mcz\\
0&bcz\
\end{pmatrix}\in Z(K).
$$ 
Otherwise, it follows from \textbf{2)} that the relation $0\neq (mcz)z'\in N$  is true for some $z'\in Z$. Then
$$
0\neq \begin{pmatrix}
a&m\\
0&b\\
\end{pmatrix}\begin{pmatrix}
czz'&0\\
0&czz'\\
\end{pmatrix} = \begin{pmatrix}
aczz'&(mcz)z'\\
0&bczz'\\
\end{pmatrix}\in Z(K).
$$ 
It remains to consider the case $a = 0$, $b\neq 0$. By condition \textbf{1)}, we have $0\neq bc\in \text{Ann}_IM\cap Z$ for some $c\in Z$. If $mc\in N$, then
$$
0\neq \begin{pmatrix}
0&m\\
0&b\\
\end{pmatrix}\begin{pmatrix}
c&0\\
0&c\\
\end{pmatrix} = \begin{pmatrix}
0&mc\\
0&bc\
\end{pmatrix}\in Z(K).
$$ 
Otherwise, it follows from \textbf{2)} that $0\neq (mc)z\in N$ for some $z\in Z$. Then 
$$
0\neq \begin{pmatrix}
0&m\\
0&b\\
\end{pmatrix}\begin{pmatrix}
cz&0\\
0&cz\\
\end{pmatrix} = \begin{pmatrix}
0&mcz\\
0&bcz\
\end{pmatrix}\in Z(K).
$$ 
Consequently, $K$ is a centrally essential ring.~$\square$

\section{Proof of Theorem 1.7}

\textbf{3.1. Lemma \cite[Lemma 3.3]{KHHC19}.} Let $R$ be a ring and let $s, t\in Z(R)$. Then
$$
Z(H_{(s,t)}(R)) = \left\lbrace\left.\begin{pmatrix}
a&0&0\\
c&d&e\\
0&0&f
\end{pmatrix}\right.\in H_{(s,t)}(R)\; \mid \; c, e, f\in Z(R)\right\rbrace.
$$
 \textbf{3.2. Lemma.} If $H_{(s,t)}(R)$ is a centrally essential ring, then $R$ is a centrally essential ring.

\textbf{Proof.} Let $0\neq r\in R$. Since $H_{(s,t)}(R)$ is a centrally essential ring, we have that for a matrix 
$A = \begin{pmatrix}
r&0&0\\
0&r&0\\
0&0&r
\end{pmatrix}$ (in this case, we have $c = e = 0$ and $a = d = f = r$), there exists a matrix 
$B = \begin{pmatrix}
a'&0&0\\
c'&d'&e'\\
0&0&f'
\end{pmatrix}\in Z(H_{(s,t)}(R))$ such that
$$
0\neq AB = \begin{pmatrix}
r&0&0\\
0&r&0\\
0&0&r
\end{pmatrix}\begin{pmatrix}
a'&0&0\\
c'&d'&e'\\
0&0&f'
\end{pmatrix} = \begin{pmatrix}
ra'&0&0\\
rc'&rd'&re'\\
0&0&rf'
\end{pmatrix}\in Z(H_{(s,t)}(R)),
$$
where $rc', re', rf'\in Z(R)$ and at least one of these products is not equal to zero. Otherwise, $ rd' = r(f' + te') = rf' + tre' = 0$,  $ra' = r(d' + sc') = rd' + src' = 0$. Then $AB = 0$. This is a contradiction. Consequently, $R$ is a centrally essential ring.~$\square$

\textbf{3.3. Lemma.} Let $R$ be a centrally essential ring and let $0\neq A = \begin{pmatrix}
a&0&0\\
c&d&e\\
0&0&f
\end{pmatrix}\in H_{(s,t)}(R)$. In each of the following cases, the ring $H_{(s,t)}(R)$ contains a  central matrix $B$ such that 
$0\neq AB\in Z(H_{(s,t)}(R))$:
\begin{enumerate}
\item[\textbf{1)}]
$a\neq 0$;
\item[\textbf{2)}]
$c = 0$;
\item[\textbf{3)}]
$f\neq 0$.
\end{enumerate}
\textbf{Proof.} \textbf{1)} Let $a\neq 0$. Then $d\neq -sc$. In the matrix $B = \begin{pmatrix}
a'&0&0\\
c'&d'&e'\\
0&0&f'
\end{pmatrix}$, we set $e' = f' = 0$. Then $d' = 0$ and $a' = sc'$. We have
$$
AB = \begin{pmatrix}
a&0&0\\
c&d&e\\
0&0&f
\end{pmatrix}\begin{pmatrix}
sc'&0&0\\
c'&0&0\\
0&0&0
\end{pmatrix} = \begin{pmatrix}
sac'&0&0\\
(cs + d)c'&0&0\\
0&0&0
\end{pmatrix}.
$$
Then $0\neq (cs + d)c'\in Z(R)$ for a suitable $c'\in Z(R)$ and $0\neq AB\in Z(H_{(s,t)}(R))$, $B\in Z(H_{(s,t)}(R))$.

\textbf{2)} Let $c = 0$. If $d\neq 0$, then $a\neq 0$, and we are in case \textbf{1)}. Let $d = -sc = 0$. Then $a = 0$, $f = -te$ and $e\neq 0$ (otherwise $A = 0$). In this case, we have 
$$
AB = \begin{pmatrix}
0&0&0\\
0&0&e\\
0&0&-te
\end{pmatrix}\begin{pmatrix}
a'&0&0\\
c'&d'&e'\\
0&0&f'
\end{pmatrix} = \begin{pmatrix}
0&0&0\\
0&0&ef'\\
0&0&-tef'
\end{pmatrix},
$$
where $c', e'\in Z(R)$ and $0\neq ef'\in Z(R)$ for a suitable $f'\in Z(R)$.

\textbf{3)} Assume that \textbf{1)} and \textbf{2)} do not hold. It follows from $a = 0$ that $-sc = f + te$ and $f = -sc - te\neq 0$. We set $c' = d' = 0$. Then $a' = 0$ and $f' = -te'$. We have:
$$
AB = \begin{pmatrix}
0&0&0\\
c&d&e\\
0&0&f
\end{pmatrix}\begin{pmatrix}
0&0&0\\
0&0&e'\\
0&0&f'
\end{pmatrix} = \begin{pmatrix}
0&0&0\\
0&0&de'+ef'\\
0&0&ff'
\end{pmatrix},
$$
where $ff' = -tfe'$, $de' + ef' = -sce' - tee' = (-sc - te)e' = fe'$. Then $0\neq fe'\in Z(R)$ for a suitable $e'\in Z(R)$.~$\square$

\textbf{3.4. Lemma.} If $R$ is a centrally essential ring, then for any $a, b\in R$ not simultaneously equal to zero, there exists an element $z\in Z(R)$ such that $az, bz\in Z(R)$ and $az\neq 0$ or $bz\neq 0$. In particular, for any non-zero-divisors $a, b\in R$, there exists $z\in Z(R)$ such that 
$0\neq az, bz\in Z(R)$.

\textbf{Proof.} It is sufficient to consider the case $a, b\notin Z(R)$, $a\neq 0$, $b\neq 0$. Since $R$ is centrally essential, we have that  $0\neq az\in Z(R)$ and $0\neq (az - b)z'\in Z(R)$ for some $z, z'\in Z(R)$. Then $bz'\in Z(R)$. If $azz'\neq 0$, then $bzz'\in Z(R)$. If $azz' = 0$, then $bz'\neq 0$. Now, if $az'\neq 0$, then $0\neq(az')z''\in Z(R)$ for some $z''\in Z(R)$ and $bz'z''\in Z(R)$.~$\square$

\textbf{3.5. Completion of the proof of Theorem 1.7.}

Let $H_{(s,t)}(R)$ be a centrally essential ring. Then $R$ is a centrally essential ring by Lemma 3.2.

Conversely, let $A = \begin{pmatrix}
a&0&0\\
c&d&e\\
0&0&f
\end{pmatrix}\in H_{(s,t)}(R)$. Let us find such a matrix $B = \begin{pmatrix}
a'&0&0\\
c'&d'&e'\\
0&0&f'
\end{pmatrix}\in Z(H_{(s,t)}(R))$ that $0\neq AB\in Z(H_{(s,t)}(R))$. It follows from Lemma 3.3 that such a matrix $B$ exists provided $a\neq 0$, $f\neq 0$ and 
$c = 0$. Therefore, we assume that $a = f = 0$ and $c\neq 0$.
We have:
$$
AB = \begin{pmatrix}
0&0&0\\
c&d&e\\
0&0&0
\end{pmatrix}\begin{pmatrix}
a'&0&0\\
c'&d'&e'\\
0&0&f'
\end{pmatrix} = \begin{pmatrix}
0&0&0\\
ca'+dc'&dd'&de'+ef'\\
0&0&0
\end{pmatrix},
$$
where 
$$
ca' + dc' = c(d' + sc') - scc' = cd' = cf' + tce',
$$
$$
dd' = -scd' = -sc(f' + te') = -scf' - stce',
$$
$$
de' + ef' = -sce' + ef'.
$$
We assume that $sc\neq 0$ or $tc\neq 0$. We set $f' = c' = 0$. Then $a' = d' = te'$ and $cd' = tce'$, $dd' = -stce'$, $de' = -sce'$. By Lemma 3.4, we have $(sc)e', (tc)e'\in Z(R)$ for a suitable $e'\in Z(R)$, and at least one of these products is not equal to zero. In both cases, $dd'\in Z(R)$. Consequently, 
$0\neq AB \in Z(H_{(s,t)}(R))$, where $B\in Z(H_{(s,t)}(R))$.

If $sc = 0$ and $tc = 0$, then $ca' + dc' = cf'$, $dd' = 0$, $de' + ef' = ef'$. It follows from Lemma 3.4 that for some $f'\in Z(R)$, we have that $cf', ef'\in Z(R)$ and $cf'\neq 0$ or $ef'\neq 0$. Consequently, $H_{(s,t)}(R)$ is a centrally essential ring.~$\square$

\section{Remarks and Examples}

For a  ring $R$, we denote by $C(R)$ the set of non-zero-divisors of $R$.

\textbf{4.1. Remark.} Let $R$ be a ring and let $I$ be a central ideal such that $I\cap C(R)\neq 0$. Then $R$ is a centrally essential ring.

Indeed, for $0\neq r\in R$, we have $0\neq rz\in I\subseteq Z(R)$ for every $z\in I\cap C(R)$.

\textbf{4.2. Remark.} Let $R$ be a semiprime ring. The ring $R$ is commutative if and only if $I\cap C(R)\neq 0$ for some central ideal $I$. In particular, a domain is commutative if and only if it contains a non-zero central ideal.

Indeed, if $R$ is commutative, then we may take the ring $R$ as $I$. In the converse implication, it follows from Remark 4.1 that $R$ is a centrally essential ring. Then the ring $R$ is commutative by \cite[Proposition 2.8]{MT19b}. 

\textbf{4.3. Example.} \cite[Example 5.1]{LT25IJ} We consider a non-commutative ring
$$
K = \left\lbrace\left.\begin{pmatrix}
a&m\\
0&b
\end{pmatrix}\;\;\right|\;\; a, b\in \mathbb{Z},\, a - b\in 2\mathbb{Z},\, m\in \mathbb{Z}_4 \right\rbrace.
$$
It follows from Corollary 1.4 that the ring $K$ is centrally essential. Moreover, it is shown in \cite{LT25IJ} that $K$ is a centrally essential ring such that all quotient rings of $K$ are centrally essential (such rings are called \textbf{completely centrally essential rings} in \cite{LT25IJ}).

\textbf{4.4. Remark.} Let $R$ be a ring and let $I$ be a central ideal of $R$ such that the quotient ring $R/I$ is commutative (such rings are called \textbf{$CIFC$ rings} in \cite{Jin21}). The ring 
$$
K = \left\lbrace\left.\begin{pmatrix}
a&m\\
0&b\\
\end{pmatrix}\right. \Big| \; a - b\in I;\, m\in R/I\right\rbrace
$$
is centrally essential. Indeed, with the use of the notation of Theorem 1.2, we have $N = R/I$ and $J = \text{Ann}_RR/I = I$. Therefore, the conditions 1) and 2) of the theorem are fulfilled. Consequently, $K$ is a centrally essential ring.

\textbf{4.5. Example.} We give an example of a centrally essential ring $K$ such that the ring $R$ is not centrally essential. We consider the $CIFC$ ring from \cite[Example 1.2]{Jin21}  as the ring $R$. Let $F$ be a field and let $A = \langle x, y\rangle$ be the free algebra generated by non-commuting variables $x, y$ over $F$. Let $B = \{f\in A | \, f(0, 0) = 0\}$. We consider the ideal $I$ of $A$ generated by all products $abc$, where $a, b, c\in B$. We set $R = A/I$. We identify elements of $A$ with their images in $R$ and we remark that $B^3 = 0$. Furthemore, let $J$ be the ideal of $R$ generated by $xy$ and $yx$. It is clear that $J\subseteq Z(R)$. In addition, it is proved in \cite[Example 1.2]{Jin21} that the quotient ring $R/J$ is commutative. Consequently, $R$ is a $CIFC$ ring. In addition, $R$ is not a centrally essential ring. Indeed, for example, for the element $x$, there does not exist  $z\in Z(R)$ such that $0\neq xz\in Z(R)$. However, the ring $K$
$$
K = \left\lbrace\left.\begin{pmatrix}
a&m\\
0&b\\
\end{pmatrix}\right. \Big| \; a, b\in R;\ a - b\in J;\, m\in R/J\right\rbrace
$$
is centrally essential by Remark 4.4.

\textbf{4.6. Remark.} In \cite[Example 1.7]{MT19b}, it is proved that if a ring $R$ has a central ideal $I$ such that $R/I$ is a field, then $R$ is a centrally essential ring. Example 4.5 shows that this assertion is not true provided $R/I$ is a commutative ring which is not a field.

\textbf{4.7. Example.} If a ring $R$ is centrally essential, then the ring $K$ is not necessarily centrally essential. Indeed, the ring
$$
K = \left\lbrace\left.\begin{pmatrix}
a&m\\
0&b\\
\end{pmatrix}\right. \Big| \; a, b\in \mathbb{Z}; \, a - b\in 2\mathbb{Z};\, m\in \mathbb{Z}\right\rbrace
$$
is not centrally essential by Corollary 1.4, since $\text{Ann}_IM = \text{Ann}_{2\mathbb{Z}}\mathbb{Z} = 0$.

\textbf{4.8. Example.} Let $R$ be a commutative uniserial artinian ring with Jacobson radical $J = J(R)$ of nilpotence index $n$. Then $R$ is a local ring with composition series $0 = J^n \subset J^{n-1} \subset \ldots \subset J \subset R$. We consider the ring
$$
K = \left\lbrace\left.\begin{pmatrix}
a&m\\
0&b\\
\end{pmatrix}\right | \; a, b\in R;\, m\in J \,\, \mbox{and} \,\, a - b\in J \right\rbrace.
$$
We have $\text{Ann}_IM = \text{Ann}_MI = J^{n-1}\le_e J$. By Corollary 1.3, the ring $K$ is centrally essential. 
In particular, specializing $R$, we obtain a sequence of finite non-commutative centrally essential subrings $K$ of $T(R, M)$ for for each prime integer $p$:
$$
K = \left\lbrace\left.\begin{pmatrix}
a&m\\
0&b\\
\end{pmatrix}\right | \; a, b\in \mathbb{Z}_{p^k};\, m\in p\mathbb{Z}_{p^k} \,\, \mbox{and} \,\, a - b\in p\mathbb{Z}_{p^k}\right\rbrace.
$$
The ring $K$ is non-commutative for $k > 2$. The smallest ring of the specified form has $128$ elements for $p = 2$ and $k = 3$.


\begin{thebibliography}{99}

\bibitem{BPR12} Birkenmeier G.F., Park J.K., Rizvi S.T. // Generalized Triangular Matrix Rings and the Fully Invariant Extending Property // Rocky Mountain Journal of Mathematics. -- 2012. -- Vol. 32, no. 4. -- P.~1299--1319.

\bibitem {DJHM25} Danchev P., Javan A., Hasanzadeh O., Moussavi A. Rings in which all elements are the sum of a central element and an element from 
$\Delta(R)$ // arXiv:2503.03643 [math.RA], https://doi.org/10.48550/arXiv.2503.03643.

\bibitem{HV99} Haghany A., Varadarajan K. // Study of formal triangular matrix rings // Comm. Algebra. -- 1999. -- Vol. 27, no. 11. -- 
P.~5507--5525.

\bibitem{Jin21} Jin H., Kim N., Lee Y., Piao Z., Ziembowski M. Structure of rings with commutative factor rings for some ideals contained in their centers
// Hacet. J. Math. Stat. -- 2021. -- Vol. 50, no. 5. -- P.~1280--1291.

\bibitem{KHHC19} Kurtulmaz Y., Halicioglu S., Harmanci A., Chen H // Rings in which elements are a sum of a central and a unit element // Bull. Belg. Math. Soc. Simon Stevin. -- 2019. -- Vol. 26, no. 4. -- P.~619--631.

\bibitem{LT25IJ} Lyubimtsev O.V., Tuganbaev A.A. Completely Centrally Essential Rings // International Journal of Algebra and Computation. -- 2025. --
Vol.~35, no 6. -- P.~927--937.

\bibitem{MT19b} Markov V.T., Tuganbaev A.A. Rings essential over their centers // Comm. Algebra. -- 2019. -- Vol. 47, no. 4. -- P.~189--194.

\bibitem{S86} Sakano K. Classical quotient rings of trivial extensions // Tsukuba J. Math. -- 1986. -- Vol. 10, no. 1. -- P.~175--181.

\end{thebibliography}
\end{document}